\documentclass[12pt, a4paper]{article}
\usepackage[pagewise]{lineno}
\usepackage{graphicx}
\usepackage{amssymb}
\usepackage{latexsym, bm}
\usepackage{multicol}
\usepackage{indentfirst}
\usepackage{amssymb,amsfonts}
\usepackage{amsmath}
\usepackage{cases}
\usepackage[noadjust]{cite}
\usepackage[colorlinks,citecolor=red]{hyperref}
\newtheorem{theorem}{Theorem}
\newtheorem{lemma}{Lemma}

\newtheorem{proposition}{Proposition}
\newtheorem{definition}{Definition}

\newtheorem{remark}{Remark}

\usepackage{amssymb}
\numberwithin{equation}{section}
\numberwithin{theorem}{section}
\numberwithin{remark}{section}
\numberwithin{definition}{section}
\numberwithin{lemma}{section}
\numberwithin{corollary}{section}
\numberwithin{proposition}{section}
\numberwithin{notation}{section}
\title{Liouville theorems for the fractional Navier-Stokes equations with arbitrary asymptotic state at infinity}
\author{Changzhi Liu\quad  Wenke Tan\footnote{tanwenkeybfq@163.com(W.Tan) ~liuchangzhi97@163.com(C.Liu)}\\
{\small Key Laboratory of Computing and Stochastic Mathematics (Ministry of Education),}\\
{\small School of Mathematics and Statistics, Hunan Normal University,}\\
{\small Changsha, Hunan 410081, China}\\
}

\date{}

\begin{document}
\maketitle
{\bf Abstract:} We mainly consider a Liouville-type problem for the three dimensional stationary fractional Navier-Stokes equations with arbitrary asymptotic state $u_\infty$ at infinity. When $u_\infty\neq 0$ and $\frac{1}{2}\leq s<1$, we prove a complete Liouville theorem by establishing some refined $L^p$ estimates for the velocity without relying on perturbation arguments. These new estimates are stronger than the $L^3$ estimates obtained by the classical perturbation framework, we thus can take $u$ as a test function and give a direct and simple proof of Liouville theorem while avoiding some technical fractional calculus. When $u_\infty\neq 0, s=\frac{1}{2}$ or $u_\infty=0,\frac{1}{2}\leq s\leq\frac{5}{6}$, we also prove a complete Liouville theorem by using frequency localization to overcome the obstacles coming from the non-local effects of $(-\Delta)^s$.
We wish to emphasize that our method dealing with the case of $u_\infty=0$ is also applicable to dimension $n$ with $n\geq 2$ and $\frac{1}{2}\leq s\leq \frac{n+2}{6}$.

\medskip

\section{Introduction}
In this paper, we consider the the stationary fractional Navier-Stokes equations
\begin{equation}\label{FNS}
 \left\{\begin{array}{ll}
(-\Delta)^s u+\nabla\cdot(u\otimes u)+\nabla P=0\quad \quad s>0,\\
\nabla\cdot u=0.
\end{array}\right.
\end{equation} in the whole space $\mathbb{R}^3$, where the unknowns $u(x)=(u_1(x),u_2(x),u_3(x))$ and $P(x)$ denote the velocity vector and the
scalar pressure at the point $x=(x_1,x_2,x_3)\in\mathbb{R}^3$, respectively. For the non-local operator $(-\Delta)^s$, it is defined for any $f\in\mathcal{S}'(\mathbb{R}^3)$ through the
Fourier transform: $(-\Delta)^sf(x)=\mathcal{F}^{-1}(|\xi|^{2s}\mathcal{F}f(\xi))$.
 Generally, we deal with solutions $u$ of \eqref{FNS} in the class of the finite Dirichlet integral
\begin{align}\label{D}
D(u)=\int_{\mathbb{R}^3}|(-\Delta)^\frac{s}{2}u|^2dx<\infty
\end{align} with the uniform condition at infinity
\begin{align}\label{far}
\lim_{|x|\to\infty}u(x)=u_\infty.
\end{align}
It is well known that the system \eqref{FNS} enjoys the scaling symmetry $(u, P)\mapsto (u^\lambda,P^\lambda)$ for any $\lambda>0$, where
\begin{align}\label{symmetry}
u^\lambda(x)=\lambda^{2s-1}u(\lambda x);\quad\quad P^\lambda(x)=\lambda^{4s-2}P(\lambda x)
\end{align}Among other things, this means that when $\frac{1}{2}\leq s<\frac{5}{6}$, the Dirichlet energy \eqref{D} is supercritical, when $s=\frac{5}{6}$ the Dirichlet energy \eqref{D} is critical, and the Dirichlet energy \eqref{D} is subcritical when $s>\frac{5}{6}$.

when $s=1$, the system \eqref{FNS} becomes the classical stationary Navier-Stokes equations.
\begin{equation}\label{NS}
 \left\{\begin{array}{ll}
-\Delta u+\nabla\cdot(u\otimes u)+\nabla P=0,\\
\nabla\cdot u=0.
\end{array}\right.
\end{equation}
Due to the pioneering work of Leray \cite{Leray}, there is a famous Liouville-type problem that whether $u\equiv0$ is the
only smooth solution of \eqref{NS} obeying the conditions \eqref{D} (called D-solution) and \eqref{far} for $u_\infty=0$. This interesting question is explicitly discussed in Galdi's book (Remark
X.9.4, \cite{Galdi}; see also Conjecture 2.5 in Tsai's book \cite{Tsai1}).

This Liouville-type problem,
without any additional assumptions, remains widely open. One of the attempts made to solve the above or related problems was firstly presented by Galdi in \cite{Galdi}, where he showed that if $u\in L^{\frac{9}{2}}(\mathbb{R}^3)$, then it holds that $u\equiv0$.
Chae and Wolf \cite{C-W} gave a logarithmic improvement of Galdi's
result by assuming that
\begin{align*}
\int_{\mathbb{R}^3}|u|^\frac{9}{2}\{\ln(2+|u|^{-1})\}^{-1}dx<\infty.
\end{align*} Also, Chae \cite{C} showed that the condition
\begin{align*}
\Delta u\in L^\frac{6}{5}(\mathbb{R}^3)
\end{align*} is sufficient for $u=0$ in $\mathbb{R}^3$. He emphasizes that the norm $\Delta u$ in $L^\frac{6}{5}(\mathbb{R}^3)$ corresponds to that of $\nabla u$ in $L^2(\mathbb{R}^3)$ at the level of scaling and that there is no mutual implication relation between their results \cite{C} and \cite{Galdi}. Kozono, Terasawa and Wakasugi proved in \cite{Kozono} that $u=0$ if the vorticity satisfies
\begin{align*}
\limsup_{|x|\to \infty}|x|^\frac{5}{3}|\omega(x)|\leq(\delta D(u))^\frac{1}{3}
\end{align*} or the velocity satisfies
\begin{align*}
||u||_{L^{\frac{9}{2},\infty}(\mathbb{R}^3)}\leq (\delta D(u))^\frac{1}{3}
\end{align*} for a small constant $\delta$. Then, the restriction imposed on the norm $||u||_{L^{\frac{9}{2},\infty}(\mathbb{R}^3)}$ in \cite{Kozono} was relaxed by Seregin and Wang in \cite{S-W}. For other values of the parameter
$p\neq\frac{9}{2}$ in the Lebesgue spaces $L^p(\mathbb{R}^3)$,  Chamorro, Jarr\'in and Lemari\'e-Rieusset \cite{C-J-L} prove that $u\in L^p(\mathbb{R}^3)$ with $3\leq p<\frac{9}{2}$
implies $u=0$. When $\frac{9}{2}< p<6$, the authors in \cite{C-J-L} also showed that $u=0$ if $u\in L^p(\mathbb{R}^3)\cap \dot{B}^{\frac{3}{p}-\frac{3}{2}}_{\infty,\infty}(\mathbb{R}^3)$, where $\dot{B}^{\frac{3}{p}-\frac{3}{2}}_{\infty,\infty}(\mathbb{R}^3)$ is a homogeneous Besov space. Moreover, for the value $p=6$, an interesting result of Seregin given
in \cite{S} shows that this problem is solved in the space $u\in BMO^{-1}(\mathbb{R}^3)\cap L^6(\mathbb{R}^3)$. Jarr\'{i}n \cite{J} established Liouville type theorem for $u\in L^p(\mathbb{R}^3)$ with $\frac{3}{2}<p<3$ and $u\in L^p(\mathbb{R}^3)\cap \dot{H}^{-1}(\mathbb{R}^3)$ with $\frac{9}{2}<p<\infty$. 

Recently, some results with some local conditions were established by several authors. Seregin \cite{S2} proved that $u=0$ if
\begin{align*}
\sup_{R>0}R^\beta(\frac{1}{|B(R)|}\int_{B(R)}|u|^qdx)^\frac{1}{q}<\infty
\end{align*} for $\frac{3}{2}<q<3$ and $\beta>\frac{6q-3}{8q-6}$. Moreover, he also arrived at the same conclusion under the
assumption
\begin{align*}
\sup_{R>0}R^{\gamma-\frac{3}{s}}||u||_{L^{s,\infty}(B(R))}<\infty,
\end{align*}where $2<s\leq3$ and $\gamma>\frac{4q-3}{6q-6}$.
 Seregin and Wang \cite{S-W} proved the vanishing of $u$ assuming either for $3<q<\infty$, $3\leq l\leq\infty$ (or $q=l=3$),
 \begin{align*}
 \liminf_{R\to\infty}R^{\frac{2}{3}-\frac{2}{q}}||u||_{L^{q,l}(B(R)\setminus B(\frac{R}{2}))}\leq \delta D(u)^\frac{3}{2},
 \end{align*}with $\delta$ a small constant, or for $\frac{12}{5}<q<3$, $1\leq l\leq\infty$, $\gamma>\frac{1}{3}+\frac{1}{q}$,
\begin{align*}
 \liminf_{R\to\infty}R^{\gamma-\frac{2}{q}}||u||_{L^{q,l}(B(R)\setminus B(\frac{R}{2}))}<\infty.
 \end{align*} Tsai \cite{Tsai} established the Liouville type theorem, using the assumption
\begin{align*}
 \liminf_{R\to\infty}\frac{1}{R}||u||^{3-\delta}_{L^q(R\leq|x|\leq\lambda R)}=0,
 \end{align*} where $q=q(\delta)=\frac{3-\delta}{1-\frac{\delta}{6}}$, for some constants $0\leq\delta\leq1$ and $\lambda>1$, or the assumption
\begin{align*}
 \liminf_{R\to\infty}R^\beta||u||^{3-\delta}_{L^q(R\leq|x|\leq R+R^{1-\alpha})}=0,
 \end{align*}where $q=q(\delta)=\frac{3-\delta}{1-\frac{\delta}{6}}$ and  $\beta=\max\{\frac{(3-\alpha)q^{-1}-2+3\alpha}{2-\delta}, \frac{-1+2\alpha}{3-\delta}\}$, for some constants $0\leq\delta\leq1$ and $\alpha\geq0$. Cho, Neustupa and Yang \cite{C-N-Y} proved that if $u$ satisfies
 \begin{align*}
 \liminf_{R\to\infty}R^{\frac{1}{3}-\frac{2}{q}}||u||_{L^q(B(2R)\setminus B(R))}<\infty\quad for \quad some\quad \frac{3}{2}<q<3
 \end{align*} or
 \begin{align*}
 \liminf_{R\to\infty}R^{\frac{-1}{3}}||u||_{L^3(B(2R)\setminus B(R))}=0,
 \end{align*}then $u=0$. For more information on the Liouville-type theorems for the classical stationary Navier-Stokes equations \eqref{NS}, one can refer to \cite{Chae1,Chae,C-Y,C-W1,CY,Yuan} and the references therein.

Compared with the classical stationary Navier-Stokes equations, there are very few results on the steady fractional Navier-Stokes equations \eqref{FNS}. Under some mild assumptions over the external force, there exists at least one solution $u\in \dot{H}^s(\mathbb{R}^3)$ obtained by Chamorro-Poggi in \cite{Ch-Po}.
Tang and Yu \cite{TY} studied partial H\"{o}lder regularity of the steady fractional Navier-Stokes equations with $\frac{1}{2}<s<1$, and proved that if $\frac{1}{2}<s<\frac{5}{6}$,
 any steady suitable weak solution is regular away from a relatively closed set with zero $5-6s$-Hausdorff measure and it is regular when $\frac{5}{6}\leq s<1$.
Corresponding to the classical Navier-Stokes equations, there is a natural Liouville-type problem for system \eqref{FNS}: When $0<s<1$, whether the smooth solution $u$ to \eqref{FNS} satisfying
the conditions \eqref{D} and \eqref{far} must be $u_\infty$.
 When $u_\infty=0$, by using the Caffarelli-Sivestre extension \cite{Caff} to convert the non-local operator $(-\Delta)^s$
on $\mathbb{R}^3$ to the local operator $\bar{\Delta}$ on $\mathbb{R}^4_+$, Wang-Xiao \cite{W-X} proved that smooth solutions to \eqref{FNS}
vanish identically when $u\in \dot{H}^s(\mathbb{R}^3)\cap L^\frac{9}{2}(\mathbb{R}^3)$, and later Yang \cite{Yang} extended this result for range $\frac{5}{6}\leq s<1$. We also refer to the recent $L^p-$theory results by Jarr\'in-Vergara-Hermosilla in \cite{J-V-H} and the corresponding result of \cite{Yang} in high dimension by Liu-Zuo \cite{L-Z}. On the other
hand, for $\frac{3}{10}<s<1$, Chamorro-Poggi \cite{Ch-Po} addressed the Liouville-type problem by assuming $u\in\dot{H}^s(\mathbb{R}^3)\cap L^{p(s)}(\mathbb{R}^3)$, where the parameter $p(s)$ depends on $s$ and approaches the critical value $\frac{6}{3-2s}$ in a certain sense. Recently, Tan \cite{Tan} obtained a Liouville type theorem for $\frac{1}{2}<s\leq 1$ by assuming $u\in \dot{H}^s(\mathbb{R}^3)\cap \dot{B}^{1-2s}_{\infty,\infty}(\mathbb{R}^3)$.

It is worth pointing out that if we consider the uniform condition \eqref{far} quipped with a non-zero state $u_\infty\neq0$, the integrability condition of $u$ as $L^\frac{9}{2}(\mathbb{R}^3)$ or $\dot{B}^{1-2s}_{\infty,\infty}(\mathbb{R}^3)$ is unreasonable. Recently, by using the strategy of proving the Liouville theorem for the case $s=1$ in \cite{Galdi}, Wang-Yang-Yu \cite{W-Y-Y} investigated this issue with $u_\infty=(1,0,0)^T$ and prove that for $\frac{1}{2}<s<1$, the smooth solution $u\in\dot{H}^s(\mathbb{R}^3)\cap \dot{W}^{1+2s,\infty}(\mathbb{R}^3)$ must equal to $u_\infty$. For more results on conditional Liouville properties for $u_\infty\neq0$, we refer to \cite{Finn,La} and the references therein.

In summary, the Liouville-type problem for fractional Navier-Stokes equations \eqref{FNS}, without any additional assumptions as $u\in L^\frac{9}{2}(\mathbb{R}^3)$ or $u\in \dot{W}^{1+2s,\infty}(\mathbb{R}^3)$, remains widely open.
In this paper, our main object is to solve the above Liouville-type problem with arbitrary asymptotic state $u_\infty$ at infinity.

Our main theorem reads:
\begin{theorem}\label{main}
Assume that $u$ be a smooth solution of \eqref{FNS} in the class \eqref{D} with \eqref{far}. It holds that
\\
(1) if $\frac{1}{2}\leq s<1$ and $u_\infty\neq0$, then $u\equiv u_\infty$,\\
(2) if $\frac{1}{2}\leq s\leq\frac{5}{6}$ and $u_\infty=0$, then $u\equiv0$.
\end{theorem}
\begin{remark}Some remarks are in order\\
(1) When $s=1$ and $u_\infty\neq 0$, the Liouville-type problem was solved by Galdi \cite{Galdi}. In this theorem, we completely solve the Liouville-type problem for $\frac{1}{2}\leq s<1$.\\
(2) When $u_\infty=0$, the Liouville-type problem becomes more complicated, we completely solve the Liouville-type problem for $\frac{1}{2}\leq s\leq\frac{5}{6}$. Our results are also applicable to the $n$-dimensional fractional Navier-Stokes equations for $n\geq 2$ with $\frac{1}{2}\leq s\leq\frac{n+2}{6}$. When $n\geq 4$, we can choose $s=1$ and our results recover the Liouville theorem for the classical stationary Navier-Stokes equations \cite{Galdi}. We wish to emphasize that in most existing works, the authors localized the nonlocal operator $(-\Delta)^s$ by Caffarelli–Silvestre extension \cite{Caff} for $s<1$ and by Young extension \cite{YR} for $s>1$. In our paper, when $n\geq 5$, we establish a unified framework to deal with the nonlocal operator $(-\Delta)^s$ for $\frac{1}{2}\leq s<1$ and $1<s\leq \frac{n+2}{6}$.\\
(3) After preparation of the manuscript, Lee and Lee posted a preprint \cite{Lee} in which they prove a similar result to (2) in Theorem \ref{main} with $\frac{n}{6}\leq s<\frac{n+2}{6}$ and $3\leq n\leq 6$. We wish to emphasize that we analyze the Liouville-type problem in the frequency spaces, but they analyze this problem in physical spaces. Our results are also applicable to $\mathbb{R}^n$ for $n\geq 2$. Furthermore, our scope of $s$ is $\frac{1}{2}\leq s\leq\frac{n+2}{6}$ which includes the critical case $s=\frac{n+2}{6}$ and the left endpoint of our scope is $\frac{1}{2}$ which is independent of $n$ and strict less than $\frac{n}{6}$ when $n\geq 4$.  
\end{remark}

To prove Theorem \ref{main}, we establish some new estimates for the solution to system \eqref{FNS}. Our new estimates which are stronger than the $L^3$ estimate established in \cite{Galdi} and \cite{W-Y-Y} give a direct and simple proof of Liouville theorem while avoiding some technical fractional calculus.
\begin{proposition}\label{kl}
Assume $u$ be a smooth solution of \eqref{FNS} in the class \eqref{D} with \eqref{far} and $u_\infty\neq0$. Denote $v=u-u_\infty$, we have the follows.\\
(1) When $\frac{5}{6}<s<1$, denote $K=[\log_4\frac{1}{6s-5}]+1$, $s_{K}=4^{K}(s-\frac{5}{6})+\frac{5}{6}$ and $\theta=\frac{1-s}{s_K-s}$, it holds
\begin{align}\label{pr1}
||v||_{\dot{H}^1(\mathbb{R}^3)}\leq C||v||^{2^K\theta+(1-\theta)}_{\dot{H}^s(\mathbb{R}^3)}.
\end{align}
(2) When $\frac{1}{2}<s\leq\frac{5}{6}$, denote $L=[\log_3\frac{1}{2s-1}]+1$, $s_L=3^L(s-\frac{1}{2})+\frac{1}{2}$ and $\theta=\frac{1-s}{s_L-s}$, it holds
\begin{align}\label{pr2}
||v||_{\dot{H}^1(\mathbb{R}^3)}\leq C||v||^{\theta L}_{L^\infty(\mathbb{R}^3)}||v||_{\dot{H}^s(\mathbb{R}^3)}.
\end{align}
(3) Assume $\frac{1}{2}<s<1$. If $\frac{2s(1+s)}{3(1-s)}\notin\mathbb{Z}$, denote $N=[\frac{2s(1+s)}{3(1-s)}]+1$, it holds
\begin{align}\label{pr3}
||v||_{L^{r}(\mathbb{R}^3)}\leq C||v||^{N}_{\dot{H}^1(\mathbb{R}^3)}||v||_{\dot{H}^s(\mathbb{R}^3)}
\end{align}where $r<2$ is defined by $\frac{1}{r}=\frac{3-2s}{6}+N(\frac{1}{2}-\frac{s}{1+s}).$ If $\frac{2s(1+s)}{3(1-s)}\in\mathbb{Z}$, denote $N=\frac{2s(1+s)}{3(1-s)},$ it follows
\begin{align}\label{pr4}
||v||_{L^{2}(\mathbb{R}^3)}\leq C||v||^{N}_{\dot{H}^1(\mathbb{R}^3)}||v||_{\dot{H}^s(\mathbb{R}^3)}.\end{align}
\end{proposition}
\begin{remark}Some remarks are in order\\
(1) To the best of our knowledge, the estimates \eqref{pr3} and \eqref{pr4} are completely new. An interesting observation in Proposition \ref{kl} is that when we want to obtain a lower integrability of $v$, we should firstly lift the regularity of $v$. \\
(2) It is worth pointing out that our estimates \eqref{pr3} and \eqref{pr4} are not applicable to the cases where $s=\frac{1}{2},u_\infty\neq0$ or $u_\infty=0$, the estimates \eqref{pr3} and \eqref{pr4} emerge from the interplay of the nonlocal operator $(-\Delta)^s$ and the nonzero asymptotic behavior $u_\infty$. The invalidity of \eqref{pr3} and \eqref{pr4} is also the reason why we must adopt some new methods to deal with the case of $u_\infty=0$ or $s=\frac{1}{2}$.
\end{remark}

We now outline the proof of Theorem \ref{main} and highlight some new ideas involved. To solve the Liouville-type problem for stationary fractional Navier-Stokes equations, the crucial and challenging step is to establish a lower integrability for solutions. When $u_\infty\neq 0$ and $\frac{1}{2}<s<1$, denoting $v=u-u_\infty$, it follows from the system \eqref{FNS} and the conditions \eqref{D}, \eqref{far} that
\begin{equation}
 \left\{\begin{array}{ll}
(-\Delta)^s v+v\cdot\nabla v+\partial_{x_1} v+\nabla P=0,\notag\\
\nabla\cdot v=0,\notag
\end{array}\right.
\end{equation}Observing the equations of $v$, we know that the nonlinear term $v\cdot\nabla v$ plays a role reducing the integrability, but the fractional linear Oseen system plays a role lifting the integrability. Precisely, assume $v\in L^r(\mathbb{R}^3)$ and $\nabla v\in L^p(\mathbb{R}^3)$, the nonlinear term $v\cdot\nabla v$ reduces the integrability such that $v\cdot\nabla v\in L^q(\mathbb{R}^3)$ with $\frac{1}{q}=\frac{1}{r}+\frac{1}{p}$. The loss of integrability is measured by $r-q>0.$ Then, the fractional linear Oseen system lifts the integrability such that $v\in L^{r_1}(\mathbb{R}^3)$ with $\frac{1}{r_1}=\frac{1}{q}-\frac{s}{1+s}$. The gain of integrability is measured by $r_1-q>0$. Noting that
\begin{align*}
\frac{1}{r_1}=\frac{1}{q}-\frac{s}{1+s}=\frac{1}{r}+(\frac{1}{p}-\frac{s}{1+s}),
\end{align*}when $p=\frac{1+s}{s}$, it follows that $r=r_1$, the influence causing by the nonlinear term is balanced by the regularization effect of the fractional linear Oseen system. If $\nabla v\in L^\frac{1+s}{s}$, the authors in \cite{W-Y-Y} adopt a delicate perturbation argument which has been used to deal with the classical Navier-Stokes equations in \cite{Galdi} to establish the $L^3$ estimate of $v$. Unlike the approach in \cite{W-Y-Y}, our strategy is based on a more direct and simple observation that if we choose $p<\frac{s+1}{s}$, we directly obtain a loss of integrability then an iterative process may deduce the required lower order integrability of $v$. Note that $\frac{1+s}{s}>2$ for $s<1$, we can choose $p=2$ and the remaining task is to establish, by exploiting certain subcritical feature of solutions, the $\dot{H}^1$ estimate of $v$. Based on the obtained lower estimate, we can choose $v$ as a test function and give a simple proof of the Liouville theorem. 

When $s=\frac{1}{2}$, since $\dot{H}^\frac{1}{2}(\mathbb{R}^3)$ is embedding into $L^3(\mathbb{R}^3)$, following the classical framework of Galdi \cite{Galdi} (also see \cite{W-Y-Y}), one can choose $v(x)\varphi_R(x)$ as a test function where $\varphi_R(x)$ is a standard cut-off function which belongs to $C^\infty_0(B_{2R})$ and equals $1$ in $B_R$. In this case, the nonlinear term $v\cdot\nabla v$ is not a trouble, the main obstacle comes from the non-local operator $(-\Delta)^s$. More precisely, the obstacle comes from the following commutator estimates:
\begin{align*}
\int_{\mathbb{R}^3}(-\Delta v)^s(v\varphi_R)dx=\int_{\mathbb{R}^3}|(-\Delta)^\frac{s}{2}v|^2\varphi_Rdx+O(R^{-(s-\frac{1}{2})}),
\end{align*}when $s=\frac{1}{2}$, the remainder $O(R^{-(s-\frac{1}{2})})$ is not vanishing as $R\to\infty$. To overcome this obstacle, we observe that the high-frequency part of velocity $v^k$, whose definition can be found in Definition \ref{LPD}, is suitable for a test function. By this choice, it is natural to consider the Liouville-type problem in the frequency spaces. We thus apply the frequency localization to analyze the different contribution comes from the different frequency part in the Dirichlet energy. Using this strategy, we establish the Liouville theorems for the case of $s=\frac{1}{2}$, $u_\infty\neq 0$ and the case of $\frac{1}{2}\leq s\leq\frac{5}{6}$, $u_\infty=0$.

Our paper is organized as follows: In Section 2, we present some notations and recall some necessary
Lemmas. In Section 3, we present the proofs of the main theorem.

\section{Preliminaries}

{\bf Notation:} In this paper, we denote $B_r=\{x\in\mathbb{R}^3:|x|\leq r\}$.  The matrix $\nabla u$ denotes the gradient of $u$ with respect to the $x$ variable, whose
$(i,j)$-th component is given by $(\nabla u)_{i,j}=\partial_ju_i$ with $1\leq i,j\leq3$. Throughout this
paper, $C$ stands for some positive constant, which may vary from line to line. Given a Banach space $X$, we denote its norm by $||\cdot||_{X}$.
We use $\mathcal{S}(\mathbb{R})^3$ and $\mathcal{S}'(\mathbb{R}^3)$ to denote Schwartz functions and the tempered distributions spaces on $\mathbb{R}^3$, respectively. Denote $L^p(\mathbb{R}^3)$ by the usual Lebesgue space with the norm
\begin{equation}
||f||_{L^p(\mathbb{R}^3)} \left\{\begin{array}{ll}
(\int_{\mathbb{R}^3}|f|^pdx)^\frac{1}{p},\quad 1\leq p<\infty,\notag\\
\text{ess} \sup_{x\in\mathbb{R}^3}|f(x)|,\quad p=\infty.\notag
\end{array}\right.
\end{equation}

Next, we will recall some well-known facts about the Littlewood-Paley decomposition. Firstly, let us recall that for every $f\in\mathcal{S}(\mathbb{R}^3)$
the Fourier transform of $f$ is defined by
\begin{align*}
(\mathcal{F}f)(\xi)=\hat{f}(\xi)=(2\pi)^{-\frac{3}{2}}\int_{\mathbb{R}^3}e^{-ix\cdot\xi}f(x)dx\quad for\quad\xi\in\mathbb{R}^3
\end{align*}The inverse Fourier transform of $g\in\mathcal{S}(\mathbb{R}^3)$ is given by
\begin{align*}
(\mathcal{F}^{-1}g)(x)=\breve{g}(x)=(2\pi)^{-\frac{3}{2}}\int_{\mathbb{R}^3}e^{ix\cdot\xi}g(\xi)dx\quad for \quad x\in\mathbb{R}^3.
\end{align*}By using Fourier transform, we can define homogeneous Sobolev norm $||\cdot||_{\dot{H}^s(\mathbb{R}^3)}$ with $s\in\mathbb{R}$ as
\begin{align*}
||f||_{\dot{H}^s(\mathbb{R}^3)}=(\int_{\mathbb{R}^3}|\xi|^{2s}|\hat{f}(\xi)|^2d\xi)^\frac{1}{2}.
\end{align*}

We now introduce the Littlewood-Paley decomposition of distributions into dyadic blocks of frequencies:
\begin{definition}\label{LPD}
Let $\psi(\xi)\in C^\infty_0(B_1)$ be a non-negative function so that $\psi(\xi)=1$ for $|\xi|\leq \frac{1}{2}$. Let $\varphi(\xi)$ be defined as $\varphi(\xi)=\psi(2^{-1}\xi)-\psi(\xi)$.
For given $u\in\mathcal{S}'(\mathbb{R}^3)$, the homogeneous dyadic blocks $\dot{\Delta}_k$ and the homogeneous low-frequency cut-off operator $\dot{S}_k$ are defined for all $k\in\mathbb{Z}$ by
\begin{align}\label{block}
\dot{\Delta}_ku(x)&=\mathcal{F}^{-1}(\varphi(2^{-k}\cdot)\hat{u}(\cdot))(x)=\frac{1}{(2\pi)^\frac{3}{2}}\int_{\mathbb{R}^3}e^{ix\cdot\xi}\varphi(2^{-k}\xi)\hat{u}(\xi)d\xi,\\
\dot{S}_ku(x)&=\mathcal{F}^{-1}(\psi(2^{-k}\cdot)\hat{u}(\cdot))(x)=\frac{1}{(2\pi)^\frac{3}{2}}\int_{\mathbb{R}^3}e^{ix\cdot\xi}\psi(2^{-k}\xi)\hat{u}(\xi)d\xi.
\end{align}
\end{definition}

 Throughout this paper, we will use the notation that $\tilde{\dot{\Delta}}_lu=\sum_{|l'-l|\leq2}\dot{\Delta}_{l'}u$. In the homogeneous case, the following Littlewood-Paley decomposition makes sense
\begin{align*}
u(x)=\sum_{k\in\mathbb{Z}}\dot{\Delta}_ku(x)\quad for\quad u\in\mathcal{S}_h'(\mathbb{R}^3),
\end{align*} where $\mathcal{S}_h'(\mathbb{R}^3)$ is given by
\begin{align*}
\mathcal{S}_h'(\mathbb{R}^3)=\{u\in\mathcal{S}'(\mathbb{R}^3):\lim_{k\to-\infty}\dot{S}_ku=0\}.
\end{align*}Moreover, it holds that
\begin{align*}
\dot{S}_ku(x)=\sum_{l\leq k-1}\dot{\Delta}_lu(x).
\end{align*}

A great advantage of the localized techniques in frequency is the so-called Bernstein inequalities which will be used in the sequel.
\begin{lemma}\cite{B-C-D}
Let $\mathcal{C}$ be an annulus and $B$ a ball. A constant C exists such that
for any nonnegative integer $k$, any couple $(p,q)$ in $[1,\infty]^2$ with $1\leq p\leq q$, and any function $u\in L^p$, we have
\begin{align}
supp\hat{u}\subset\lambda B&\Rightarrow ||\nabla^k u||_{L^q}\leq C^{k+1}\lambda^{k+\frac{3}{p}-\frac{3}{q}}||u||_{L^p},\label{Ber1}\\
supp\hat{u}\subset\lambda \mathcal{C}&\Rightarrow C^{-k-1}\lambda^k||u||_{L^p}\leq||\nabla^ku||_{L^p}\leq C^{k+1}\lambda^k||u||_{L^p}\label{Ber2}
\end{align}
\end{lemma}
Based on the Littlewood-Paley decomposition and Bernstein inequality \eqref{Ber2}, the homogeneous Sobolev norm $||\cdot||_{\dot{H}^s}$ can be equivalently written as
\begin{align*}
||f||_{\dot{H}^s(\mathbb{R}^3)}=(\sum_{k\in\mathbb{Z}}2^{2ks}||\dot{\Delta}_k f||^2_{L^2})^\frac{1}{2}.
\end{align*}

When we investigate the nonlinear term in the stationary fractional Navier-Stokes equations, we need the following estimate.
\begin{lemma}\label{bds1}\cite{B-C-D}
Assume $0<s<1$ and $f,g\in\dot{H}^s(\mathbb{R}^3)$. Then it holds that
\begin{align}
||f g||_{\dot{H}^{2s-\frac{3}{2}}(\mathbb{R}^3)}\leq C||f||_{\dot{H}^s(\mathbb{R}^3)}||g||_{\dot{H}^s(\mathbb{R}^3)}.
\end{align}
\end{lemma}
It is also worth pointing out that throughout this paper, we will use the following well-known facts:
\begin{align}\label{+10}
||\dot{S}_kf||_{L^p}\leq C||f||_{L^p},\quad ||f^k||_{L^p}\leq C||f||_{L^p}\quad for \quad 1\leq p\leq\infty
\end{align}here $f^k=f-\dot{S}_kf$.

Finally, we introduce the following multiplier theorem by Lizorkin \cite{Liz}.
\begin{theorem}\label{P}(\text{Lizorkin}).
Let $\Phi: \mathbb{R}^n\rightarrow\mathbb{C}$ be continuous together with the derivative $$\frac{\partial^n\Phi}{\partial{\xi_1}...\partial{\xi_n}}$$ and all preceding derivatives for
$|\xi_i|>0$, $i=1,...,n$. Then, if for some $\beta\in[0,1)$ and $M>0$
\begin{align}\label{P1}
|\xi_1|^{\kappa_1+\beta}\cdot...\cdot|\xi_n|^{\kappa_n+\beta} \left| \frac{\partial^\kappa\Phi}{\partial\xi^{\kappa_1}_1...\partial\xi^{\kappa_n}_n} \right| \leq M,
\end{align}where $\kappa_i$ is zero or one and $\kappa=\sum_{i=1}^{n}\kappa_i$, the integral transform
\begin{align*}
Tu(x)=\frac{1}{(2\pi)^\frac{n}{2}}\int_{\mathbb{R}^n}e^{ix\cdot\xi}\Phi(\xi)\hat{u}(\xi)d\xi,\quad u\in\mathcal{S}(\mathbb{R}^n),
\end{align*}
defines a bounded linear operator from $L^q(\mathbb{R}^n)$ to $L^r(\mathbb{R}^n)$, $1<q<\infty$, $\frac{1}{r}=\frac{1}{q}-\beta$, and we have
$$||Tu||_{L^r(\mathbb{R}^n)}\leq C(q,\beta,M)||u||_{L^q(\mathbb{R}^n)}.$$
\end{theorem}

\section{Proofs of Theorem \ref{main}}
We divide the proof of Theorem \ref{main} into four steps.\\
{\bf Step I: the proof of Proposition \ref{kl}.} Firstly, from the orthogonal transform invariant of the system \eqref{FNS}, we may assume $u_\infty=(C,0,0)^T$ with $C\neq0$. Then, by the following scaling invariant of fractional Navier-Stokes system
 \begin{align*}
 u_\lambda(x)=\lambda^{2s-1}u(\lambda x),\quad P_\lambda(x)=\lambda^{4s-2}P(\lambda x)
 \end{align*} we can assume $u_\infty=(1,0,0)^T$.

 Denoting $v=u-u_\infty$, it follows from the system \eqref{FNS} and the conditions \eqref{D}, \eqref{far} that
\begin{equation}\label{MFNS}
 \left\{\begin{array}{ll}
(-\Delta)^s v+\nabla\cdot(v\otimes v)+\partial_{x_1} v+\nabla P=0,\\
\nabla\cdot v=0,
\end{array}\right.
\end{equation} and
\begin{align}\label{b1+}
||v||_{\dot{H}^s(\mathbb{R}^3)}<\infty,
\end{align}
\begin{align}\label{far1}
\lim_{|x|\to\infty}v=0.
\end{align}

Next, we focus on the generalized fractional Ossen system \eqref{MFNS}. We will firstly derive the $\dot{H}^1$ estimate for velocity. For the case of  $\frac{5}{6}<s<1$, applying the Fourier transform, we obtain from \eqref{MFNS} that
\begin{equation}\label{fourier}
 \left\{\begin{array}{ll}
(|\xi|^{2s} +i\xi_1)\hat{v}+i\xi\hat{P}=-i\xi\cdot(\widehat{v\otimes v}),\\
i\xi\cdot\hat{v}=0.
\end{array}\right.
\end{equation} Taking $L^2$ inner product to $\eqref{fourier}_1$ with $i\xi$, we get
\begin{align}\label{fp}
\hat{P}=-\frac{\xi\cdot(\widehat{v\otimes v})\cdot\xi}{|\xi|^2}.
\end{align}Substituting \eqref{fp} into \eqref{fourier}, one can see that
\begin{align}\label{fv}
\hat{v}_j(\xi)=-(|\xi|^{2s}+i\xi_1)^{-1}[(\delta_{jk}-\frac{\xi_j\xi_k}{|\xi|^2})i\xi_l \widehat{v_k v_l} ].
\end{align}

Multiplying \eqref{fv} by $|\xi|^\sigma$, we get
\begin{align*}
|\xi|^\sigma\hat{v}_j(\xi)=-\frac{|\xi|^{\sigma-2s+\frac{3}{2}}i\xi_l}{|\xi|^{2s}+i\xi_1}(\delta_{jk}-\frac{\xi_j\xi_k}{|\xi|^2})|\xi|^{2s-\frac{3}{2}}\widehat{v_k v_l}.
\end{align*}By taking $\sigma=4s-\frac{5}{2}$, one can see that
\begin{align*}
|\xi|^\sigma\hat{v}_j=\Phi_{jkl}(\xi)|\xi|^{2s-\frac{3}{2}}\widehat{v_k v_l},
\end{align*}with $\Phi_{jkl}(\xi)=-\frac{|\xi|^{2s-1}i\xi_l}{|\xi|^{2s}+i\xi_1}(\delta_{jk}-\frac{\xi_j\xi_k}{|\xi|^2})$. Direct calculation shows that the Fourier multiplier $\Phi_{jkl}(\xi)$ satisfies the conditions in Theorem \ref{P} with $\beta=0$. Indeed, let $$\Phi_{jkl}(\xi):=\Phi_1+\Phi_2$$with
$\Phi_1=-\frac{|\xi|^{2s-1}i\xi_l}{|\xi|^{2s}+i\xi_1}\delta_{jk}$ and $\Phi_2=\frac{|\xi|^{2s-1}i\xi_l}{|\xi|^{2s}+i\xi_1}\frac{\xi_j\xi_k}{|\xi|^2}$,
it suffices to consider $\Phi_1$, since the calculations for $\Phi_2$ are similar. It is clear that
\begin{align*}
|\Phi_1|\leq \frac{|\xi|^{2s}}{\sqrt{|\xi|^{4s}+|\xi_1|^2}}\leq 1,
\end{align*} and
\begin{equation*}
    \begin{aligned}
        \partial_{\xi_h} \Phi_{1} = -\left[ \frac{(2s-1) \lvert \xi \rvert^{2s-2} i \xi_l \frac{ \xi_h }{\lvert \xi \rvert} + \lvert \xi \rvert^{2s-1} i \delta_{hl}}{\lvert \xi \rvert^{2s} + i \xi_1} - \frac{(\lvert \xi \rvert^{2s-1} i \xi_l) (2s \lvert \xi \rvert^{2s-1} \frac{ \xi_h }{\lvert \xi \rvert} + i \delta_{h1} ) }{(\lvert \xi \rvert^{2s} + i \xi_1)^2}   \right]      \delta_{jk}.
    \end{aligned}
\end{equation*}
Furthermore, we get by using Young inequality that
\begin{equation*}
    \begin{aligned}
        &\left| \xi_h\right|  \left|\partial_{\xi_h} \Phi_{1}  \right|\\ = &\left| \frac{(2s-1) \lvert \xi \rvert^{2s-3} i \xi_l \xi_h \left| \xi_h \right|  + \lvert \xi \rvert^{2s-1} \left| \xi_h\right| i \delta_{hl}}{\lvert \xi \rvert^{2s} + i \xi_1} - \frac{(2s \lvert \xi \rvert^{4s-3} i \xi_l \xi_h \left| \xi_h \right| -  \lvert \xi \rvert^{2s-1} \xi_l \lvert \xi_h \rvert i \delta_{h1} ) }{(\lvert \xi \rvert^{2s} + i\xi_1)^2}   \right|      \delta_{jk}\\
        \leq &2(\frac{|\xi|^{2s}}{\sqrt{|\xi|^{4s}+|\xi_1|^2}}+\frac{|\xi|^{4s}}{|\xi|^{4s}+|\xi_1|^2}+\frac{|\xi|^{2s}|\xi_1|}{|\xi|^{4s}+|\xi_1|^2}) \\
        \leq& 6.   \end{aligned}
\end{equation*} This proves \eqref{P1} when $\kappa=\kappa_1+\kappa_2+\kappa_3=1$, the proofs for $\kappa=2$ or $\kappa=3$ are similar. Combining this together, we obtain that
\begin{align*}
|\xi_1|^{\kappa_1}  |\xi_2|^{\kappa_2}  |\xi_3|^{\kappa_3} \left| \frac{\partial^\kappa\Phi_{jkl}}{\partial\xi^{\kappa_1}_1  \partial\xi^{\kappa_2}_2  \partial\xi^{\kappa_3}_3 } \right|\leq M.
\end{align*}
for all $\kappa=\sum_{j=1}^3\kappa_j$  and $\kappa_j\in\{0,1\}$.

 Taking $r = q = 2$ and $\beta = 0$ in Theorem \ref{P}, we have that
\begin{align}\label{siagm}
||v||_{\dot{H}^\sigma(\mathbb{R}^3)}\leq C||v\otimes v||_{\dot{H}^{2s-\frac{3}{2}}(\mathbb{R}^3)}\leq C||v||^2_{\dot{H}^s(\mathbb{R}^3)},
\end{align}where we have used Lemma \ref{bds1}. Since $\frac{5}{6}<s<1$, it follows that $$\sigma-s=3s-\frac{5}{2}>0.$$  We thus obtain a gain of regularity.

Next, from this gain of regularity, we can construct an iteration scheme as follows.  Denote $s_1=\sigma$. If $s_1\geq1$, we terminate this scheme. Otherwise, noting that $\sigma=4s-\frac{5}{2}<\frac{3}{2}$, we replace $s$ by $s_1$ and repeat the above arguments to get $s_2$ satisfying $$s_2=4s_1-\frac{5}{2}=4^2s-4\frac{5}{2}-\frac{5}{2}$$ and
\begin{align*}
||v||_{\dot{H}^{s_2}}\leq C||v||^2_{\dot{H}^{s_1}}\leq C^3||v||^4_{\dot{H}^s(\mathbb{R}^3)}
 \end{align*}where we have used \eqref{siagm}.
 Repeating this iteration scheme finite times, we finally obtain $K=[\log_4\frac{1}{6s-5}]+1$ and $s_{K}=4^{K}(s-\frac{5}{6})+\frac{5}{6}$ such that $s_K\geq1$ and $v$ obeying
 \begin{align*}
 ||v||_{\dot{H}^{s_K}(\mathbb{R}^3)}\leq C^{2^K-1}||v||^{2^{K}}_{\dot{H}^s(\mathbb{R}^3)}.
 \end{align*} From this fact and \eqref{b1+}, we deduce by using interpolation that
\begin{align}\label{h1}
||v||_{\dot{H}^1(\mathbb{R}^3)}\leq C^{\theta(2^K-1)}||v||^{2^K\theta+(1-\theta)}_{\dot{H}^s(\mathbb{R}^3)}
\end{align}here $\theta=\frac{1-s}{s_K-s}$. We thus get \eqref{pr1}.

Now we focus on the case of $\frac{1}{2}<s\leq\frac{5}{6}$. Observing from the smoothness of $v$ and the uniform condition \eqref{far1} that
\begin{align}\label{es1+}
||v||_{L^\infty}<\infty.
\end{align} From \eqref{es1+} and \eqref{b1+}, we have that
\begin{align*}
v\in \dot{H}^s(\mathbb{R}^3)\cap L^\infty(\mathbb{R}^3).
\end{align*}  Since $\dot{H}^s(\mathbb{R}^3)\cap L^\infty(\mathbb{R}^3)$ is an algebra for $s>0$, we conclude that
\begin{align}\label{alg}
||v\otimes v||_{\dot{H}^s(\mathbb{R}^3)}\leq C||v||_{L^\infty(\mathbb{R}^3)}||v||_{\dot{H}^s(\mathbb{R}^3)}.
\end{align} By using the similar arguments as the case of $s\in(\frac{5}{6},1)$, one can get that
\begin{align*}
|\xi|^\delta\hat{v}_j(\xi)=-\frac{|\xi|^{\delta-s}i\xi_l}{|\xi|^{2s}+i\xi_1}(\delta_{jk}-\frac{\xi_j\xi_k}{|\xi|^2})|\xi|^{s}\widehat{v_k v_l}.
\end{align*}By taking $\delta=3s-1$, it follows that
\begin{align*}
|\xi|^\delta \hat{v}_j(\xi)=\Psi_{jkl}(\xi)|\xi|^s\widehat{v_k v_l},
\end{align*}with $\Psi_{jkl}(\xi)=-\frac{|\xi|^{2s-1}i\xi_l}{|\xi|^{2s}+i\xi_1}(\delta_{jk}-\frac{\xi_j\xi_k}{|\xi|^2})$. By applying Theorem \ref{P} with $\beta=0, r=q=2$, we also obtain that
\begin{align}\label{delta}
||v||_{\dot{H}^\delta(\mathbb{R}^3)}\leq C||v\otimes v||_{\dot{H}^s(\mathbb{R}^3)}\leq C||v||_{L^\infty(\mathbb{R}^3)}||v||_{\dot{H}^s(\mathbb{R}^3)}
\end{align}where we have used \eqref{alg}. Since $\frac{1}{2}<s\leq \frac{5}{6}$, it holds that
$$\delta-s=2s-1>0.$$ We also get a gain of regularity. By using a similar arguments as the case of $\frac{5}{6}<s<1$, we also get that $L=[\log_3\frac{1}{2s-1}]+1$ and $s_L=3^L(s-\frac{1}{2})+\frac{1}{2}$ such that $s_L\geq 1$ and $v$ obeying
\begin{align*}
||v||_{\dot{H}^{s_L}(\mathbb{R}^3)}\leq C^L||v||^L_{L^\infty(\mathbb{R}^3)}||v||_{\dot{H}^s(\mathbb{R}^3)}.
\end{align*}By using interpolation, we have
\begin{align}\label{h12}
||v||_{\dot{H}^1(\mathbb{R}^3)}\leq C^{\theta L}||v||^{\theta L}_{L^\infty(\mathbb{R}^3)}||v||_{\dot{H}^s(\mathbb{R}^3)}
\end{align} here $\theta=\frac{1-s}{s_L-s}$. We thus obtain \eqref{pr2}.

With the $\dot{H}^1$ estimate \eqref{h1}(or \eqref{h12}) in hand, we can deduce the low order integrability of $v$. Firstly, we rewrite the system \eqref{MFNS} as
\begin{equation}\label{LMFNS}
 \left\{\begin{array}{ll}
(-\Delta)^s v+\partial_{x_1} v+\nabla P=-v\cdot\nabla v,\\
\nabla\cdot v=0.
\end{array}\right.
\end{equation} The linear part of the equations
of $v$ are fractional Oseen equations but fractional Stokes equations. This fact is crucial for our analysis. From Sobolev embedding theorem $\dot{H}^s(\mathbb{R}^3)\hookrightarrow L^{\frac{6}{3-2s}}(\mathbb{R}^3)$ and \eqref{b1+}, we have
\begin{align}\label{emb}
||v||_{L^{\frac{6}{3-2s}}(\mathbb{R}^3)}\leq C||v||_{\dot{H}^s(\mathbb{R}^3)}<\infty.
\end{align} From \eqref{h1}(or \eqref{h12}), \eqref{emb} and  H\"{o}lder inequlity, we see that
\begin{align}\label{nli}
||v\cdot\nabla v||_{L^{\frac{3}{3-s}}(\mathbb{R}^3)}\leq ||v||_{L^{\frac{6}{3-2s}}(\mathbb{R}^3)}||v||_{\dot{H}^1(\mathbb{R}^3)}<\infty.
\end{align}

By performing the  Fourier transform on the system \eqref{MFNS}, we can solve $v$ from \eqref{LMFNS}:
\begin{align}\label{sol}
\hat{v}(\xi)=\frac{1}{|\xi|^{2s}+i\xi_1}(Id-\frac{\xi\otimes\xi}{|\xi|^2})\hat{f}=\Gamma(\xi)\hat{f}
\end{align}with $\Gamma(\xi)=\frac{1}{|\xi|^{2s}+i\xi_1}(Id-\frac{\xi\otimes\xi}{|\xi|^2})$ and $f=-v\cdot\nabla v$. By using a similar calculation dealing with $\Phi_1$ (also see Lemma 2.7 in \cite{W-Y-Y}),
one can obtain that
\begin{align}\label{ga}
|\xi_1|^{\kappa_1+\beta}|\xi_2|^{\kappa_2+\beta}|\xi_3|^{\kappa_3+\beta} \left| \frac{\partial^\kappa\Gamma(\xi)}{\partial^{\kappa_1}\xi_1\partial^{\kappa_2}\xi_2\partial^{\kappa_3}\xi_3} \right| \leq M
\end{align}with $\beta=\frac{s}{1+s}$ and $\kappa=\sum_{j=1}^{3}\kappa_j$, $\kappa_j\in\{0,1\}$.
Since $s\in(\frac{1}{2},1)$ implies $1<\frac{3}{3-s}<\frac{1+s}{s}$, we can apply Theorem \ref{P} to the formula \eqref{sol} then deduce that
\begin{align}\label{low1}
||v||_{L^{r_1}(\mathbb{R}^3)}\leq C||v\cdot\nabla v||_{L^{\frac{3}{3-s}}(\mathbb{R}^3)}\leq C_1||v||_{L^\frac{6}{3-2s}(\mathbb{R}^3)}||v||_{\dot{H}^1(\mathbb{R}^3)}
\end{align}with $\frac{1}{r_1}=\frac{3-s}{3}-\frac{s}{1+s}=\frac{3-2s}{6}+(\frac{1}{2}-\frac{s}{1+s})$. A simple and useful observation is that when $s\in (\frac{1}{2},1)$, it holds that $\frac{1}{2}-\frac{s}{1+s}>0$ and $r_1<\frac{6}{3-2s}$. We thus get a loss of integrability. This fact inspires us that an iterative process may further lower the integrability of $v$. Indeed, combine \eqref{low1} with \eqref{h1}(or \eqref{h12}) implies that
\begin{align*}
||v\cdot\nabla v||_{L^{q_1}(\mathbb{R}^3)}\leq||v||_{L^{r_1}(\mathbb{R}^3)}||v||_{\dot{H}^1(\mathbb{R}^3)}.
\end{align*}with $\frac{1}{q_1}=\frac{1}{r_1}+\frac{1}{2}=\frac{3-s}{3}+(\frac{1}{2}-\frac{s}{1+s})>\frac{s}{1+s}$. By using Theorem \ref{P}, \eqref{emb} and \eqref{low1}, we deduce that
\begin{align}\label{low2}
||v||_{L^{r_2}(\mathbb{R}^3)}&\leq C_2||v\cdot\nabla v||_{L^{q_1}(\mathbb{R}^3)}\leq C_2||v||_{L^{r_1}(\mathbb{R}^3)}||v||_{\dot{H}^1(\mathbb{R}^3)}\\&\leq C_1C_2||v||_{\dot{H}^s(\mathbb{R}^3)}||v||^2_{\dot{H}^1(\mathbb{R}^3)}\notag
\end{align}with $\frac{1}{r_2}=\frac{1}{q_1}-\frac{s}{1+s}=\frac{3-2s}{6}+2(\frac{1}{2}-\frac{s}{1+s})$.

Iterating this process, we will get the following facts. If $\frac{2s(1+s)}{3(1-s)}\notin\mathbb{Z}$, there exists
 two sequences $\{r_k\}_{k=1}^{N+1}$ and $\{q_k\}_{k=0}^N$ obeying
\begin{align*}
N=[\frac{2s(1+s)}{3(1-s)}],~~~~~~~~~~~~~~~&\\
\frac{s}{1+s}<\frac{1}{q_k}=\frac{3-s}{3}+k(\frac{1}{2}-\frac{s}{1+s})<1,&\\
\frac{1}{r_k}=\frac{1}{q_{k-1}}-\frac{s}{1+s}=\frac{3-2s}{6}+k(\frac{1}{2}-\frac{s}{1+s})\geq&\frac{3-2s}{6}.
\end{align*}
Furthermore, we have that $v\in L^{r_{N+1}}(\mathbb{R}^3)$ and
\begin{align}\label{last}
||v||_{L^{r_{N+1}}(\mathbb{R}^3)}\leq \prod_{i=1}^{N+1}C_i||v||^{N+1}_{\dot{H}^1(\mathbb{R}^3)}||v||_{\dot{H}^s(\mathbb{R}^3)}.
\end{align}
From the formula
\begin{align*}
 \frac{1}{r_{N+1}}=\frac{3-2s}{6}+(N+1)(\frac{1}{2}-\frac{s}{1+s})>\frac{3-2s}{6}+\frac{2s(1+s)}{3(1-s)}(\frac{1}{2}-\frac{s}{1+s})=\frac{1}{2},
 \end{align*}we see that
 \begin{align}\label{index}
 r_{N+1}<2.
 \end{align} We thus get the estimate \eqref{pr3}. If $\frac{2s(1+s)}{3(1-s)}\in\mathbb{Z}$, we then obtain two sequences $\{r_k\}_{k=1}^{N}$ and $\{q_k\}_{k=0}^{N-1}$ satisfying
\begin{align*}
N=\frac{2s(1+s)}{3(1-s)},~~~~~~~~~~~~~~~&\\
\frac{s}{1+s}<\frac{1}{q_k}=\frac{3-s}{3}+k(\frac{1}{2}-\frac{s}{1+s})<1,&\\
\frac{1}{r_k}=\frac{1}{q_{k-1}}-\frac{s}{1+s}=\frac{3-2s}{6}+k(\frac{1}{2}-\frac{s}{1+s})\geq&\frac{3-2s}{6}.
\end{align*}
Furthermore, we have that $v\in L^{r_{N}}(\mathbb{R}^3)$ and
\begin{align}\label{last1}
||v||_{L^{r_{N}}(\mathbb{R}^3)}\leq \prod_{i=1}^{N}C_i||v||^{N}_{\dot{H}^1(\mathbb{R}^3)}||v||_{\dot{H}^s(\mathbb{R}^3)}.
\end{align}
From the formula
\begin{align*}
 \frac{1}{r_{N}}=\frac{3-2s}{6}+N(\frac{1}{2}-\frac{s}{1+s})=\frac{3-2s}{6}+\frac{2s(1+s)}{3(1-s)}(\frac{1}{2}-\frac{s}{1+s})=\frac{1}{2},
 \end{align*}we also see that
 \begin{align}\label{index1}
 r_{N}=2.
 \end{align}
 We thus get the estimate \eqref{pr4} then complete the proof of Proposition \ref{kl}.
 
{\bf Step II: the proof of Theorem \ref{main} for the case of $u_\infty\neq0$ and $\frac{1}{2}<s<1$.} 
From Proposition \ref{kl} and interpolation, we have that $v\in \dot{H}^1(\mathbb{R}^3)\cap L^p(\mathbb{R}^3)$ for $2\leq p\leq 6$. This means that we can take $L^2$ inner product to the system \eqref{MFNS} with $v$ and get
 \begin{align*}
 \int_{\mathbb{R}^3}|(-\Delta)^\frac{s}{2}v|^2dx=-\int_{\mathbb{R}^3}v\cdot\nabla v\cdot vdx-\int_{\mathbb{R}^3}\partial_{x_1}v\cdot vdx.
 \end{align*}By using integration by part and the fact $\nabla\cdot v=0$, one can get
 \begin{align*}
 \int_{\mathbb{R}^3}|(-\Delta)^\frac{s}{2}v|^2dx=
 \frac{1}{2}\int_{\mathbb{R}^3}|v|^2\nabla\cdot vdx-
 \frac{1}{2}\int_{\mathbb{R}^3}\partial_{x_1}|v|^2dx=0
  \end{align*}
 and thus $v=0$. Since $u=v+u_\infty$, it follows $u=u_\infty$.

{\bf Step III: the proof of Theorem \ref{main} for the case of $u_\infty\neq0$ and $s=\frac{1}{2}$.} Firstly, we recall that $v^k:=v-\dot{S}_kv$ is the high frequency part of $v$. By the definition of $\dot{H}^\frac{1}{2}$ space, we deduce that
\begin{align}\label{b4}
||\nabla v^k||_{\dot{H}^{-\frac{1}{2}}}\leq C|| v^k||_{\dot{H}^{1-\frac{1}{2}}}\leq C||v||_{\dot{H}^\frac{1}{2}}.
\end{align} On the other hand, from \eqref{alg}, we also see that 
\begin{align}\label{b3}
v\otimes v\in \dot{H}^\frac{1}{2}(\mathbb{R}^3).
\end{align}
Based on \eqref{b1+}, \eqref{b3} and \eqref{b4}, we now can take $L^2$ inner product to \eqref{MFNS} with $v^k$ then obtain that
\begin{align*}
&\int_{\mathbb{R}^3}|(-\Delta)^\frac{1}{4} v^k|^2dx\\
=&\int_{\mathbb{R}^3}v_j v_i\partial_j v^k_idx+
\int_{\mathbb{R}^3}v_i\partial_1 v^k_idx-\int_{\mathbb{R}^3}(-\Delta)^\frac{1}{4}\dot{S}_kv\cdot(-\Delta)^\frac{1}{4}v^kdx.
\end{align*}
Splitting $v=\dot{S}_k v+v^k$, we have from $\nabla\cdot v=0$ and integration by parts that
\begin{align*}
&\int_{\mathbb{R}^3}|(-\Delta)^\frac{1}{4} v^k|^2dx\\
=&-\int_{\mathbb{R}^3}\dot{S}_kv_j\partial_j \dot{S}_kv_iv^k_idx-\int_{\mathbb{R}^3}v^k_j\partial_j \dot{S}_kv_iv^k_idx+
\int_{\mathbb{R}^3}\dot{S}_kv_i\partial_1 v^k_idx\\&-\int_{\mathbb{R}^3}(-\Delta)^\frac{1}{4}\dot{S}_kv\cdot(-\Delta)^\frac{1}{4}v^kdx.
\end{align*}It is not difficult to see upon \eqref{b1+} and \eqref{b4} that
\begin{align*}
&\lim_{k\to-\infty}\int_{\mathbb{R}^3}|(-\Delta)^\frac{1}{4} v^k|^2dx=\int_{\mathbb{R}^3}|(-\Delta)^\frac{1}{4}v|^2dx,\\
&\lim_{k\to-\infty}|\int_{\mathbb{R}^3}(-\Delta)^\frac{1}{4}\dot{S}_kv\cdot(-\Delta)^\frac{1}{4}v^kdx|\leq C\lim_{k\to-\infty}||\dot{S}_k v||_{\dot{H}^\frac{1}{2}}||v||_{\dot{H}^\frac{1}{2}}=0,\\
&\lim_{k\to-\infty}|\int_{\mathbb{R}^3}\dot{S}_kv_i\partial_1v^k_idx|\leq \lim_{k\to-\infty}||\dot{S}_kv||_{\dot{H}^\frac{1}{2}}||\nabla v^k||_{\dot{H}^{-\frac{1}{2}}}\leq C\lim_{k\to-\infty}||\dot{S}_kv||_{\dot{H}^\frac{1}{2}}||v||_{\dot{H}^{\frac{1}{2}}}=0
\end{align*}Taking $k\to-\infty$, it follows that
\begin{align}\label{es2}
&\int_{\mathbb{R}^3}|(-\Delta)^\frac{1}{4} v|^2dx\\
=&-\liminf_{k\to-\infty}\{\int_{\mathbb{R}^3}\dot{S}_kv_j\partial_j \dot{S}_kv_iv^k_idx+\int_{\mathbb{R}^3}v^k_j\partial_j \dot{S}_kv_iv^k_idx\},\notag\\
=&-\liminf_{k\to-\infty}\{{I}+{J}\}.\notag
\end{align}
In the remainder of this step, our objects are to prove $\lim_{k\to\infty}\int_{\mathbb{R}^3}\dot{S}_kv_j\partial_j{S}_kv_i v^k_idx=0$ and $\lim_{k\to\infty}\int_{\mathbb{R}^3}v^k_j\partial_j{S}_kv_i v^k_idx=0$ under the conditions \eqref{b1+} and \eqref{far1}. 

To exploit the structure of the term $I$ and the bound \eqref{b1+}, we apply the Bony decomposition to the product of the low-frequency component and high frequency-component $\dot{S}_kv_j v^k_i$ and hence deduce that\begin{align}\label{BI}
I
=&\int_{\mathbb{R}^3}\partial_j\dot{S}_kv_i(\sum_{l\in\mathbb{Z}}\dot{\Delta}_l\dot{S}_kv_j\dot{S}_{l-2}v^k_i+\sum_{l\in\mathbb{Z}}\dot{\Delta}_lv^k_i\dot{S}_{l-2}\dot{S}_kv_j
+\sum_{l\in\mathbb{Z}}\dot{\Delta}_l\dot{S}_kv_j\tilde{\dot{\Delta}}_lv^k_i)dx\\
=&I_1+I_{2}+I_3.\notag
\end{align} We now aim to control $I_1-I_3$. From the definitions of $\dot{\Delta}_l$, $\dot{S}_k$ and $v^k$ we see that $\dot{\Delta}_l\dot{S}_kv_j\dot{S}_{l-2}v^k_i\neq0$ means $l\leq k$ and $l\geq k+2$. This means that $\dot{\Delta}_l\dot{S}_kv_j\dot{S}_{l-2}v^k_i=0$ for all $l\in\mathbb{Z}$ and then
\begin{align}\label{BI1}
I_1=0,
\end{align}
For $I_{2}$, we first observe that $\dot{\Delta}_lv^k_i\dot{S}_{l-2}\dot{S}_kv_j\neq0$ means $l\geq k-1$ and
\begin{align*}
supp\mathcal{F}(\dot{\Delta}_lv^k_i\dot{S}_{l-2}\dot{S}_kv_j)\subset\{\xi:2^{l-2}\leq|\xi|<\frac{9}{8}2^{l+1}\}.
\end{align*} On the other hand, we also have
\begin{align*}
supp \mathcal{F}(\partial_j\dot{S}_kv_i)\subset\{\xi:|\xi|< 2^k\}.
\end{align*}Putting all this together, we conclude that
\begin{align*}
I_{2}=&\int_{\mathbb{R}^3}\sum_{l=k-1}^{k+1}\sum_{l'=l-2}^{k-1}\partial_j\dot{\Delta}_{l'}v_i\dot{\Delta}_lv^k_i\dot{S}_{l-2}\dot{S}_{k}v_jdx.
\end{align*}
In a similar vein, from the definition of $\tilde{\dot{\Delta}}_l$, $\dot{\Delta}_l$ and $v^k$, one has\begin{align*}
I_3=&\int_{\mathbb{R}^3}\sum_{l=k-3}^{k}\dot{\Delta}_l\dot{S}_kv_j\partial_j\dot{S}_kv_i\tilde{\dot{\Delta}}_lv^k_idx.
\end{align*}

We now estimate $I_2$ and $I_3$. Using \eqref{Ber1}, \eqref{+10}, H\"{o}lder's inequality and Sobolev embedding theorem $\dot{H}^\frac{1}{2}(\mathbb{R}^3)\hookrightarrow L^3(\mathbb{R}^3)$, we get that
\begin{align*}
|I_{2}|
\leq &C\sum_{l=k-1}^{k+1}\sum_{l'=l-2}^{k-1}2^{l'}||\dot{S}_{ k}v_j||_{L^\infty}||\dot{\Delta}_{l'}v_i||_{L^2}||\dot{\Delta}_lv_i||_{L^2}\\
\leq & C||\dot{S}_{ k}v_j||_{L^\infty}\sum_{l=k-1}^{k+1}\sum_{l'=l-2}^{k-1}2^{l'\frac{1}{2}}||\dot{\Delta}_{l'}v_i||_{L^2}2^{l\frac{1}{2}}||\dot{\Delta}_lv_i||_{L^2}\notag\\
\leq & C2^{k}||\dot{S}_kv_j||_{L^3}\sum_{l=k-1}^{k+1}\sum_{l'=l-2}^{k-1}2^{l'\frac{1}{2}}||\dot{\Delta}_{l'}v_i||_{L^2}2^{l\frac{1}{2}}||\dot{\Delta}_lv_i||_{L^2},\notag\\
\leq & C2^{k}||\dot{S}_kv_j||_{\dot{H}^\frac{1}{2}}\sum_{l=k-1}^{k+1}\sum_{l'=l-2}^{k-1}2^{l'\frac{1}{2}}||\dot{\Delta}_{l'}v_i||_{L^2}2^{l\frac{1}{2}}||\dot{\Delta}_lv_i||_{L^2}\\
\leq & C2^k||v||^3_{\dot{H}^\frac{1}{2}}.
\end{align*}From \eqref{b1+}, we conclude that
\begin{align}\label{BI2}
\lim_{k\to-\infty}I_2=0.
\end{align} A similar computation deduces that
\begin{align*}
|I_3|\leq C2^{k}||\dot{S}_kv_i||_{\dot{H}^\frac{1}{2}}\sum_{l=k-3}^{k}2^{l\frac{1}{2}}||\dot{\Delta}_l v_j||_{L^2}2^{l\frac{1}{2}}||\tilde{\dot{\Delta}}_lv_i||_{L^2}\leq C2^k||v||^3_{\dot{H}^\frac{1}{2}}
\end{align*}and then
\begin{align}\label{BI3}
\lim_{k\to-\infty}I_3=0.
\end{align}Substituting \eqref{BI1}-\eqref{BI3} into \eqref{BI} implies
\begin{align}\label{BIL}
\lim_{k\to-\infty}I=0.
\end{align}

We are left with the estimation of the term $J=\int_{\mathbb{R}^3}v^k_j\partial_j{S}_kv_i v^k_idx$. To exploit the structure of the term $J$ and the bound \eqref{b1+}, we now apply the Bony decomposition to the product of the high frequency-components $v^k_j v^k_i$ and hence deduce that
\begin{align}\label{J}
J=&\int_{\mathbb{R}^3}\partial_j\dot{S}_kv_iv^k_jv^k_idx,\\
=&\int_{\mathbb{R}^3}\partial_j\dot{S}_kv_i(\sum_{l\in\mathbb{Z}}\dot{\Delta}_lv^k_j\dot{S}_{l-2}v^k_i+\sum_{l\in\mathbb{Z}}\dot{\Delta}_lv^k_i\dot{S}_{l-2}v^k_j
+\sum_{l\in\mathbb{Z}}\dot{\Delta}_lv^k_i\tilde{\dot{\Delta}}_{l}v^k_j)dx,\notag\\
=&J_1+J_2+J_3.\notag
\end{align} For $J_1$, we first observe that $\dot{\Delta}_lv^k_j\dot{S}_{l-2}v^k_i\neq0\Rightarrow l\geq k+2$ and
\begin{align*}
supp\mathcal{F}(\dot{\Delta}_lv^k_j\dot{S}_{l-2}v^k_i)\subset\{\xi:2^{l-2}\leq|\xi|<\frac{9}{8}2^{l+1}\} \quad for\quad l\geq k+2.
\end{align*}On the other hand, it is also holding that 
\begin{align*}
supp \mathcal{F}(\partial_j\dot{S}_kv_i)\subset\{\xi:|\xi|< 2^k\}.
\end{align*} Putting all this together, one can see that
\begin{align}\label{J1}
J_1=&\int_{\mathbb{R}^3}\partial_j\dot{S}_kv_i\sum_{l\in\mathbb{Z}}\dot{\Delta}_lv^k_j\dot{S}_{l-2}v^k_idx\\
=&\int_{\mathbb{R}^3}\partial_j\dot{S}_kv_i\sum_{l\geq k+2}\dot{\Delta}_lv^k_j\dot{S}_{l-2}u^k_idx\notag\\
=&\int_{\mathbb{R}^3}\mathcal{F}(\partial_j\dot{S}_kv_i)\mathcal{F}(\sum_{l\geq k+2}\dot{\Delta}_lv^k_j\dot{S}_{l-2}v^k_i)d\xi\notag\\
=&0.\notag
\end{align} Applying the same arguments to $J_2$, we also obtain
\begin{align}\label{J2}
J_2=\int_{\mathbb{R}^3}\partial_j\dot{S}_kv_i\sum_{l\in\mathbb{Z}}\dot{\Delta}_lv^k_i\dot{S}_{l-2}v^k_jdx=0.
\end{align}
For $J_3$, it is not difficult to see that $\dot{\Delta}_lv^k_i\neq0\Rightarrow l\geq k-1$, we thus get that
\begin{align*}
J_3=\int_{\mathbb{R}^3}\partial_j\dot{S}_kv_i\sum_{l\in\mathbb{Z}}\dot{\Delta}_lv^k_i\tilde{\dot{\Delta}}_{l}v^k_jdx
=\int_{\mathbb{R}^3}\partial_j\dot{S}_kv_i\sum_{l\geq k-1}\dot{\Delta}_lv^k_i\tilde{\dot{\Delta}}_{l}v^k_jdx.
\end{align*}  Using \eqref{Ber1}, \eqref{+10}, H\"{o}lder's inequality and Sobolev embedding theorem $\dot{H}^\frac{1}{2}(\mathbb{R}^3)\hookrightarrow L^3(\mathbb{R}^3)$, we get
\begin{align*}
|J_3|\leq &C 2^k||\dot{S}_kv_i||_{L^\infty}\sum_{l\geq k-1}2^{-l}2^{l\frac{1}{2}}||\dot{\Delta}_lv^k_i||_{L^2}2^{l\frac{1}{2}}||\tilde{\dot{\Delta}}_{l}v^k_j||_{L^2}\\
\leq& C2^{2k}||\dot{S}_kv_i||_{L^3}\sum_{l\geq k-1} 2^{-l}2^{l\frac{1}{2}}||\dot{\Delta}_lv||_{L^2}2^{l\frac{1}{2}}||\tilde{\dot{\Delta}}_{l}v||_{L^2}\\
\leq& C2^{k}||\dot{S}_kv_i||_{\dot{H}^\frac{1}{2}}\sum_{l\geq k-1}2^{(k-l)} 2^{l\frac{1}{2}}||\dot{\Delta}_lv||_{L^2}2^{l\frac{1}{2}}||\tilde{\dot{\Delta}}_{l}v||_{L^2}\\
\leq& C2^k||v||^3_{\dot{H}^\frac{1}{2}}.
\end{align*} From \eqref{b1+}, we deduce that
\begin{align}\label{J3}
\lim_{k\to-\infty}J_3=0.
\end{align}Substituting \eqref{J1}-\eqref{J3} into \eqref{J}, one has
\begin{align}\label{BJL}
\lim_{k\to-\infty}J=0.
\end{align}

Substituting \eqref{BIL} and \eqref{BJL} into \eqref{es2}, we obtain that
\begin{align}\label{last}
\int_{\mathbb{R}^3}|(-\Delta)^\frac{s}{2} v|^2dx=0.
\end{align}
This means that $v=0$ and then $u=u_\infty$.

{\bf Step IV: the proof of Theorem \ref{main} for the case of $u_\infty=0$ and $\frac{1}{2}\leq s\leq\frac{5}{6}$.}

We first show some estimates for $u\otimes u$ to ensure the validity of taking $L^2$ inner product to \eqref{FNS} with $\dot{S}_k u$. From Lemma \ref{bds1} one has
\begin{align}\label{es9}
||u_iu_j||_{\dot{H}^{2s-\frac{3}{2}}}\leq ||u||^2_{\dot{H}^s}=D(u).
\end{align} Secondly, by using the fact $s\leq\frac{5}{6}$, we also obtain that
\begin{align}\label{es10}
||\dot{S}_k u_i||_{\dot{H}^{\frac{5}{2}-2s}}=&(\sum_{l\leq k-1}(2^{(\frac{5}{2}-3s)l}2^{ls}||\dot{\Delta}_lu_i||_{L^2})^2)^\frac{1}{2}\leq C(k)||u||_{\dot{H}^s}=C(k)D^\frac{1}{2}(u).
\end{align} It should be point out that for general dimension $n\geq 2$, \eqref{es9} and \eqref{es10} are also valid if we replace $\frac{3}{2}$ and $\frac{5}{2}$ by $\frac{n}{2}$ and $\frac{n+2}{2}$ respectively. 

From \eqref{es9} and \eqref{es10}, one can see that the integral $\int_{\mathbb{R}^3}u_ju_i\partial_j\dot{S}_k u_idx$ is well defined. Indeed, it holds that
\begin{align}\label{es11}
&\int_{\mathbb{R}^3}u_ju_i\partial_j\dot{S}_k u_idx\\=&\int_{\mathbb{R}^3}|\xi|^{2s-\frac{3}{2}}\mathcal {F}(u_ju_i)|\xi|^{\frac{3}{2}-2s}\mathcal {F}(\partial_j\dot{S}_k u_i)d\xi\notag\\
\leq&||u_iu_j||_{\dot{H}^{2s-\frac{3}{2}}}||\partial_j\dot{S}_k u_i||_{\dot{H}^{\frac{3}{2}-2s}}\notag\\
\leq&C(k)||u_iu_j||_{\dot{H}^{2s-\frac{3}{2}}}||\dot{S}_k u_i||_{\dot{H}^{\frac{5}{2}-2s}}\notag\\
\leq& C(k)D^\frac{3}{2}(u).\notag
\end{align}

Based on \eqref{es9}, \eqref{es10} and \eqref{es11}, we now can take $L^2$ inner product to \eqref{FNS} with $\dot{S}_k u$ then obtain that
\begin{align*}
&\int_{\mathbb{R}^3}|(-\Delta)^\frac{s}{2} \dot{S}_ku|^2dx\\
=&\int_{\mathbb{R}^3}u_j u_i\partial_j \dot{S}_ku_idx
-\int_{\mathbb{R}^3}(-\Delta)^\frac{s}{2}\dot{S}_ku\cdot(-\Delta)^\frac{s}{2}u^kdx\notag\\
=&\int_{\mathbb{R}^3}\dot{S}_ku_j\partial_j{S}_ku_i u^k_idx+\int_{\mathbb{R}^3}u^k_j\partial_j{S}_ku_i u^k_idx-\int_{\mathbb{R}^3}(-\Delta)^\frac{s}{2}\dot{S}_ku(-\Delta)^\frac{s}{2}u^kdx\notag
\end{align*}where we have used the fact $\nabla\cdot u=0$ and hence $\int_{\mathbb{R}^3}u_j \dot{S}_ku_i\partial_j \dot{S}_ku_idx=0$. It is not difficult to see that
\begin{align*}
\lim_{k\to\infty}|\int_{\mathbb{R}^3}(-\Delta)^\frac{s}{2}\dot{S}_ku\cdot(-\Delta)^\frac{s}{2}u^kdx|=0,\\
\lim_{k\to\infty}\int_{\mathbb{R}^3}|(-\Delta)^\frac{s}{2} \dot{S}_ku|^2dx=\int_{\mathbb{R}^3}|(-\Delta)^\frac{s}{2}u|^2dx.
\end{align*}Taking $k\to\infty$, we conclude that
\begin{align}\label{es24}
\int_{\mathbb{R}^3}|(-\Delta)^\frac{s}{2}u|^2dx=&\liminf_{k\to\infty}\{\int_{\mathbb{R}^3}\dot{S}_ku_j\partial_j{S}_ku_i u^k_idx+\int_{\mathbb{R}^3}u^k_j\partial_j{S}_ku_i u^k_idx\}\\=&\liminf_{k\to\infty}\{K+L\}.\notag
\end{align}
In the remainder of this step, our objects are to prove $\lim_{k\to\infty}\int_{\mathbb{R}^3}\dot{S}_ku_j\partial_j{S}_ku_i u^k_idx=0$ and $\lim_{k\to\infty}\int_{\mathbb{R}^3}u^k_j\partial_j{S}_ku_i u^k_idx=0$ under the conditions \eqref{D} and \eqref{far}.

To exploit the structure of the term $K$ and the bound \eqref{D}, we apply the Bony decomposition to the product of the low-frequency component and high frequency-component $\dot{S}_ku_j u^k_i$ and hence deduce that
\begin{align}\label{bd1}
K
=&\int_{\mathbb{R}^3}\partial_j\dot{S}_ku_i(\sum_{l\in\mathbb{Z}}\dot{\Delta}_l\dot{S}_ku_j\dot{S}_{l-2}u^k_i+\sum_{l\in\mathbb{Z}}\dot{\Delta}_lu^k_i\dot{S}_{l-2}\dot{S}_ku_j
+\sum_{l\in\mathbb{Z}}\dot{\Delta}_l\dot{S}_ku_j\tilde{\dot{\Delta}}_lu^k_i)dx\\
=&K_1+K_{2}+K_3.\notag
\end{align} We now aim to control $K_1-K_3$. From the definitions of $\dot{\Delta}_l$, $\dot{S}_k$ and $u^k$ we see that $\dot{\Delta}_l\dot{S}_ku_j\dot{S}_{l-2}u^k_i\neq0$ means $l\leq k$ and $l\geq k+2$. This means that $\dot{\Delta}_l\dot{S}_ku_j\dot{S}_{l-2}u^k_i=0$ for all $l\in\mathbb{Z}$ and then
\begin{align}\label{K1}
K_{1}=\int_{\mathbb{R}^3}\partial_j\dot{S}_ku_i\sum_{l\in\mathbb{Z}}\dot{\Delta}_l\dot{S}_ku_j\dot{S}_{l-2}u^k_idx=0.
\end{align}

For $K_{2}$, we first observe that $\dot{\Delta}_lu^k_i\dot{S}_{l-2}\dot{S}_ku_j\neq0$ means $l\geq k-1$ and
\begin{align*}
supp \mathcal{F}(\dot{\Delta}_lu^k_i\dot{S}_{l-2}\dot{S}_ku_j)\subset\{\xi:2^{l-2}\leq|\xi|<\frac{9}{8}2^{l+1}\}.
\end{align*} Secondly, we also have
\begin{align*}
supp \mathcal{F}(\partial_j\dot{S}_ku_i)\subset\{\xi:|\xi|< 2^k\}.
\end{align*}Putting all this together, we conclude that
\begin{align*}
K_{2}=&
\int_{\mathbb{R}^3}\partial_j\dot{S}_ku_i\sum_{l=k-1}^{k+1}\dot{\Delta}_lu^k_i\dot{S}_{l-2}\dot{S}_{k}u_jdx\\
=&\int_{\mathbb{R}^3}\sum_{l=k-1}^{k+1}\sum_{l'=l-2}^{k-1}\partial_j\dot{\Delta}_{l'}u_i\dot{\Delta}_lu^k_i\dot{S}_{l-2}\dot{S}_{k}u_jdx.
\end{align*}
\begin{align}\label{K2+}
|K_{2}|
\leq &C\sum_{l=k-1}^{k+1}\sum_{l'=l-2}^{k-1}2^{l'}||\dot{S}_{ k}u_j||_{L^\infty}||\dot{\Delta}_{l'}u_i||_{L^2}||\dot{\Delta}_lu_i||_{L^2}\\
\leq & C2^{k(1-2s)}||\dot{S}_{ k}u_j||_{L^\infty}\sum_{l=k-1}^{k+1}\sum_{l'=l-2}^{k-1}2^{l's}||\dot{\Delta}_{l'}u_i||_{L^2}2^{ls}||\dot{\Delta}_lu_i||_{L^2}\notag\\
\leq & C2^{k(1-2s)}||u||_{L^\infty}\sum_{l=k-1}^{k+1}\sum_{l'=l-2}^{k-1}2^{l's}||\dot{\Delta}_{l'}u_i||_{L^2}2^{ls}||\dot{\Delta}_lu_i||_{L^2}.\notag
\end{align} On the one hand, one can use H\"{o}lder's inequality and \eqref{D} to see that
\begin{align*}
\lim_{k\to\infty}\sum_{l=k-1}^{k+2}\sum_{l'=l-2}^{k-1}2^{l's}||\dot{\Delta}_{l'}u_i||_{L^2}2^{ls}||\dot{\Delta}_lu_i||_{L^2}\leq\lim_{k\to\infty} \sum_{i=k-3}^{k+2}(2^{ls}||\dot{\Delta}_lu||_{L^2})^2=0,\end{align*}
on the other hand, observe from the smoothness of $u$ and the uniform condition \eqref{far} that
\begin{align}\label{es1}
||u||_{L^\infty}<\infty.
\end{align}Substituting this estimates into \eqref{K2+} and observing the fact $1-2s\leq 0$, we thus get
\begin{align}\label{K2}
\lim_{k\to\infty}K_{2}=0.
\end{align}
In a similar vein, from the definition of $\tilde{\dot{\Delta}}_l$, $\dot{\Delta}_l$ and $u^k$ one has
\begin{align*}
K_3=&\int_{\mathbb{R}^3}\sum_{l=k-3}^{k}\dot{\Delta}_l\dot{S}_ku_j\partial_j\dot{S}_ku_i\tilde{\dot{\Delta}}_lu^k_idx.
\end{align*}
From \eqref{Ber1}, \eqref{+10} and H\"{o}lder's inequality, we also obtain
\begin{align*}
|K_3|\leq C2^{k(1-2s)}||u||_{L^\infty}\sum_{l=k-3}^{k}2^{ls}||\dot{\Delta}_l u_j||_{L^2}2^{ls}||\tilde{\dot{\Delta}}_lu_i||_{L^2}.
\end{align*}By using the similar arguments for $K_2$, we thus get
\begin{align}\label{K3}
\lim_{k\to\infty}K_3=0.
\end{align}
Substituting the estimates \eqref{K1}, \eqref{K2} and \eqref{K3} into \eqref{bd1}, we have
\begin{align}\label{K}
\lim_{k\to\infty}K=0.
\end{align}

We are left with the estimation of the term $L=\int_{\mathbb{R}^3}u^k_j\partial_j{S}_ku_i u^k_idx$. To exploit the structure of the term $L$ and the bound \eqref{D}, we now apply the Bony decomposition to the product of the high frequency-components $u^k_j u^k_i$ and hence deduce that
\begin{align}\label{L}
L=&\int_{\mathbb{R}^3}\partial_j\dot{S}_ku_iu^k_ju^k_idx,\\
=&\int_{\mathbb{R}^3}\partial_j\dot{S}_ku_i(\sum_{l\in\mathbb{Z}}\dot{\Delta}_lu^k_j\dot{S}_{l-2}u^k_i+\sum_{l\in\mathbb{Z}}\dot{\Delta}_lu^k_i\dot{S}_{l-2}u^k_j
+\sum_{l\in\mathbb{Z}}\dot{\Delta}_lu^k_i\tilde{\dot{\Delta}}_{l}u^k_j)dx,\notag\\
=&L_1+L_2+L_3.\notag
\end{align}
For $L_1$, we first observe that $\dot{\Delta}_lu^k_j\dot{S}_{l-2}u^k_i\neq0\Rightarrow l\geq k+2$ and
\begin{align*}
supp\mathcal{F}(\dot{\Delta}_lu^k_j\dot{S}_{l-2}u^k_i)\subset\{\xi:2^{l-2}\leq|\xi|<\frac{9}{8}2^{l+1}\} \quad for\quad l\geq k+2.
\end{align*}Secondly, it is clear that
\begin{align*}
supp \mathcal{F}(\partial_j\dot{S}_ku_i)\subset\{\xi:|\xi|< 2^k\}.
\end{align*}Putting all this together, we have
\begin{align}\label{L1}
L_1=&\int_{\mathbb{R}^3}\partial_j\dot{S}_ku_i\sum_{l\in\mathbb{Z}}\dot{\Delta}_lu^k_j\dot{S}_{l-2}u^k_idx\\
=&\int_{\mathbb{R}^3}\partial_j\dot{S}_ku_i\sum_{l\geq k+2}\dot{\Delta}_lu^k_j\dot{S}_{l-2}u^k_idx\notag\\
=&\int_{\mathbb{R}^3}\mathcal{F}(\partial_j\dot{S}_ku_i)\mathcal{F}(\sum_{l\geq k+2}\dot{\Delta}_lu^k_j\dot{S}_{l-2}u^k_i)d\xi\notag\\
=&0.\notag
\end{align} Applying the same arguments to $L_2$, we also obtain
\begin{align}\label{L2}
L_2=\int_{\mathbb{R}^3}\partial_j\dot{S}_ku_i\sum_{l\in\mathbb{Z}}\dot{\Delta}_lu^k_i\dot{S}_{l-2}u^k_jdx=0.
\end{align}
For $L_3$, it is not difficult to see that $\dot{\Delta}_lu^k_i\neq0\Rightarrow l\geq k-1$, we thus get that
\begin{align*}
L_3=\int_{\mathbb{R}^3}\partial_j\dot{S}_ku_i\sum_{l\in\mathbb{Z}}\dot{\Delta}_lu^k_i\tilde{\dot{\Delta}}_{l}u^k_jdx
=\int_{\mathbb{R}^3}\partial_j\dot{S}_ku_i\sum_{l\geq k-1}\dot{\Delta}_lu^k_i\tilde{\dot{\Delta}}_{l}u^k_jdx.
\end{align*}From \eqref{Ber1}, \eqref{+10} and H\"{o}lder's inequality, we have
\begin{align*}
|L_3|\leq &C 2^k||\dot{S}_ku_i||_{L^\infty}\sum_{l\geq k-1}2^{-2ls}2^{ls}||\dot{\Delta}_lu^k_i||_{L^2}2^{ls}||\tilde{\dot{\Delta}}_{l}u^k_j||_{L^2}\\
\leq& C2^{k(1-2s)}||u||_{L^\infty}\sum_{l\geq k-1} 2^{ls}||\dot{\Delta}_lu||_{L^2}2^{ls}||\tilde{\dot{\Delta}}_{l}u||_{L^2}.
\end{align*}Using a similar argument for $K_2$, we obtain that
\begin{align}\label{L3}
\lim_{k\to\infty}L_3=0.
\end{align}
Substituting \eqref{L1}-\eqref{L3} into \eqref{L}, we thus get
\begin{align}\label{L+}
\lim_{k\to\infty}L=0.
\end{align}Substituting \eqref{K} and \eqref{L+} into \eqref{es24}, we have
\begin{align}
\int_{\mathbb{R}^3}|(-\Delta)^\frac{s}{2}u|^2dx=0
\end{align} and then $u=u_\infty=0$.\\

\section*{Acknowledgments}
The authors are supported by the Construct Program of the Key Discipline in Hunan Province and NSFC (Grant No. 11871209), and the Hunan Provincial NSF (No. 2022JJ10033)
\\
\noindent {\bf Conflict of interest:}
We declare that we do not have any commercial or associative interest
that represents a conflict of interest in connection with the work
submitted.

\end{document}